\documentclass[a4paper]{article}

\input{tcilatex}

\begin{document}

\begin{center}
{\Large A\ Triangle Analog to Pascal's characterizing primes }

Issam Kaddoura , Marwa Zeid

Lebanese International University

School Of Art And Science

LAAS 21th International Science Meeting ,989-992, April 17 ,2015

\bigskip
\end{center}

\bigskip \textbf{Abstract} : In the conference paper we construct analouge
to Pascal's triangle that characterize primes with additional fascinating
properties.

\textbf{Introduction }:

A primality test is an algorithm for determining whether an input number is
prime. Amongst other fields of mathematics, it is used for cryptography.
Unlike factorization, primality tests do not generally give prime factors,
only stating whether the input number is prime or not.Factorization is
thought to be a computationally difficult problem, where as primality
testing is comparatively easy (its running time is polynomial in the size of
the input). A connection between the sequence of prime numbers and the
binomial coefficients is given by a well-known characterization of the prime
numbers:

\textit{Consider the entries in the pth row of Pascal's triangle, without
the initial and final entries. }

\textit{They are all divisible by p if and only if p is a prime. }[1]

In this paper we construct an analouge to pascal's triangle that
characterize primes with additional fascinating properties .

We name this triangle \textbf{Alkarkhi's triangle }to the memory of Abu Bakr
ibn Muhammad ibn al Husayn Al-Karkhi was a 10th-century mathematician and
engineer who lived at Baghdad. His three famous surviving works are :
Al-Badi' fi'l-hisab, Al-Fakhri fi'l-jabr wa'l-muqabala, and Al-Kafi
fi'l-hisab.[2] \ \ \ \ \ \ \ \ \ \ 

\textbf{Alkarkhi's triangle} is a geometrical arrangement of numbers defined
by\ the following matrix

$K(i,j)=[a_{ij}]$ \ such that :

$a_{ij}=\left\{ 
\begin{array}{c}
\binom{i}{j}+(-1)^{j+1}\text{ \ \ \ \ \ \ \ \ \ \ \ \ \ \ \ }if\text{ \ \ \ }%
i\geq j\text{\ }>0\text{\ \ \ \ \ \ \ \ \ \ \ \ \ \ \ \ \ \ \ \ \ } \\ 
0\text{ \ \ \ \ \ \ \ \ \ \ \ \ \ \ \ \ \ \ \ \ \ \ \ }if\text{ \ \ }\ \
0<i<j\text{\ \ \ \ \ \ \ \ }%
\end{array}%
\right. $

where i is the row number and j is the column number and $\binom{i}{j}$ is a
binomial coefficient .

\bigskip \textbf{Properties of Alkarkhi's Triangle}

Alkarkhi's Triangle determines the coefficients which arise in the expansion
:

$\ $%
\begin{eqnarray*}
\sum\limits_{i=1}^{n-1}K(n,i)x^{i} &=&(1+x)^{n}-\sum_{i=0}^{n}(-x)^{i} \\
&=&\sum\limits_{i=0}^{n}\binom{n}{i}\left( x\right)
^{n-i}-\sum_{i=0}^{n}(-x)^{i} \\
&=&\left( \binom{n}{0}\left( x\right) ^{n}+\binom{n}{1}\left( x\right)
^{n-1}+\binom{n}{1}\left( x\right) ^{n-2}+...+\binom{n}{n-1}\left( x\right)
^{1}+\binom{n}{n}\left( x\right) ^{0}\right) \\
&&-\left( 1-x+x^{2}-.....+(-1)^{n}x^{n}\right) \\
&=&(\binom{n}{1}-1)\left( x\right) ^{n-1}+(\binom{n}{2}+1)\left( x\right)
^{n-2}+...+(\binom{n}{n-1}+(-1)^{n-1})\left( x\right) ^{1}
\end{eqnarray*}

\textbf{Lemma 1}:\ The enteries of the $i^{th}$ row represents the
coefficients of the above expansion of the generating function and all
enteries can be constructed from the recurrence formula :\ 
\[
K(n+1,d)+(-1)^{d}=K(n,d)+k(n,d-1) 
\]

\ 

\textbf{Proof}: \ \ \ using the identity

\[
\binom{n}{d-1}+\binom{n}{d}=\binom{n+1}{d} 
\]

and the definition 
\[
K(n,d)=\binom{n}{d}+(-1)^{d+1} 
\]

with direct subistitution

\[
K(n,d-1)+(-1)^{d-1}+K(n,d)+(-1)^{d}=K(n+1,d)+(-1)^{d} 
\]

and consequently

\[
K(n+1,d)+(-1)^{d}=K(n,d)+k(n,d-1) 
\]

\bigskip

The first few enteries of \ \textbf{Al\textbf{k}arkhi's Triangle}\ \ are
given below :

$K(n,d){\small =}\left[ 
\begin{array}{ccccccccccccc}
{\small 2} & {\small 0} & {\small 0} & {\small 0} & {\small 0} & {\small 0}
& {\small 0} & {\small 0} & {\small 0} & {\small 0} & {\small 0} & {\small 0}
& {\small 0} \\ 
{\small 3} & {\small 0} & {\small 0} & {\small 0} & {\small 0} & {\small 0}
& {\small 0} & {\small 0} & {\small 0} & {\small 0} & {\small 0} & {\small 0}
& {\small 0} \\ 
{\small 4} & {\small 2} & {\small 2} & {\small 0} & {\small 0} & {\small 0}
& {\small 0} & {\small 0} & {\small 0} & {\small 0} & {\small 0} & {\small 0}
& {\small 0} \\ 
{\small 5} & {\small 5} & {\small 5} & {\small 0} & {\small 0} & {\small 0}
& {\small 0} & {\small 0} & {\small 0} & {\small 0} & {\small 0} & {\small 0}
& {\small 0} \\ 
{\small 6} & {\small 9} & {\small 11} & {\small 4} & {\small 2} & {\small 0}
& {\small 0} & {\small 0} & {\small 0} & {\small 0} & {\small 0} & {\small 0}
& {\small 0} \\ 
{\small 7} & {\small 14} & {\small 21} & {\small 14} & {\small 7} & {\small 0%
} & {\small 0} & {\small 0} & {\small 0} & {\small 0} & {\small 0} & {\small %
0} & {\small 0} \\ 
{\small 8} & {\small 20} & {\small 36} & {\small 34} & {\small 22} & {\small %
6} & {\small 2} & {\small 0} & {\small 0} & {\small 0} & {\small 0} & 
{\small 0} & {\small 0} \\ 
{\small 9} & {\small 27} & {\small 57} & {\small 69} & {\small 57} & {\small %
27} & {\small 9} & {\small 0} & {\small 0} & {\small 0} & {\small 0} & 
{\small 0} & {\small 0} \\ 
{\small 10} & {\small 35} & {\small 85} & {\small 125} & {\small 127} & 
{\small 83} & {\small 37} & {\small 8} & {\small 2} & {\small 0} & {\small 0}
& {\small 0} & {\small 0} \\ 
{\small 11} & {\small 44} & {\small 121} & {\small 209} & {\small 253} & 
{\small 209} & {\small 121} & {\small 44} & {\small 11} & {\small 0} & 
{\small 0} & {\small 0} & {\small 0} \\ 
{\small 12} & {\small 54} & {\small 166} & {\small 329} & {\small 463} & 
{\small 461} & {\small 331} & {\small 164} & {\small 56} & {\small 10} & 
{\small 2} & {\small 0} & {\small 0} \\ 
{\small 13} & {\small 65} & {\small 221} & {\small 494} & {\small 793} & 
{\small 923} & {\small 793} & {\small 494} & {\small 221} & {\small 65} & 
{\small 13} & {\small 0} & {\small 0} \\ 
{\small \ast } & {\small \ast } & {\small \ast } & {\small \ast } & {\small %
\ast } & {\small \ast } & {\small \ast } & {\small \ast } & {\small \ast } & 
{\small \ast } & {\small \ast } & {\small \ast } & {\small \ast }%
\end{array}%
\right] $

\textbf{Alkarkhi's Triangle characterizes prime numbers}

we present the following amazing property of \textbf{Alkarkhi's Triangle \ }%
to characterize primes

Theorem 1 :\ \textit{p is prime if and only if all the entires of the (}$%
p-1)^{th}$ row in \textbf{Alkarkhi's Triangle } \textit{are zero mod p }

Proof : Obviously each entery of the \textit{(}$p-1)^{th}$ row in Alkarkhi's
Triangle has the form \ 
\[
\binom{p-1}{d}+(-1)^{d+1} 
\]%
and using :%
\[
\text{For any }0\leq d\ \leq p-1\ \ :\binom{p-1}{d}+(-1)^{d+1}\equiv 0\func{%
mod}(p)\bigskip \ \text{\ if and only if p is prime .} 
\]%
\bigskip

Notice that all the enteries in \textit{(}$p-1)^{th}$ row in bold divisible
by the number p indicates the primality of p

\[
\left[ 
\begin{array}{c}
\text{\ 2 \ Prime} \\ 
\text{ \ 3\ Prime } \\ 
\text{ 4\ \ Not\ }\ \text{Prime} \\ 
\text{ 5\ \ \ Prime} \\ 
\text{6\ Not\ Prime} \\ 
\text{\ 7 \ Prime} \\ 
\text{8 \ Not\ Prime} \\ 
\text{9 Not Prime} \\ 
\text{10 Not Prime } \\ 
\text{11 \ Prime } \\ 
\text{12 \ Not\ \ Prime } \\ 
\text{13 \ \ Prime \ } \\ 
\ast \ast \ast \ast \ast \ast \ast \ast \ast \ast%
\end{array}%
\begin{array}{ccccccccccccc}
\mathbf{2} & 0 & 0 & 0 & 0 & 0 & 0 & 0 & 0 & 0 & 0 & 0 & 0 \\ 
\mathbf{3} & 0 & 0 & 0 & 0 & 0 & 0 & 0 & 0 & 0 & 0 & 0 & 0 \\ 
4 & 2 & 2 & 0 & 0 & 0 & 0 & 0 & 0 & 0 & 0 & 0 & 0 \\ 
\mathbf{5} & \mathbf{5} & \mathbf{5} & 0 & 0 & 0 & 0 & 0 & 0 & 0 & 0 & 0 & 0
\\ 
6 & 9 & 11 & 4 & 2 & 0 & 0 & 0 & 0 & 0 & 0 & 0 & 0 \\ 
\mathbf{7} & \mathbf{14} & \mathbf{21} & \mathbf{14} & \mathbf{7} & 0 & 0 & 0
& 0 & 0 & 0 & 0 & 0 \\ 
8 & 20 & 36 & 34 & 22 & 6 & 2 & 0 & 0 & 0 & 0 & 0 & 0 \\ 
9 & 27 & 57 & 69 & 57 & 27 & 9 & 0 & 0 & 0 & 0 & 0 & 0 \\ 
10 & 35 & 85 & 125 & 127 & 83 & 37 & 8 & 2 & 0 & 0 & 0 & 0 \\ 
\boldsymbol{11} & \boldsymbol{44} & \boldsymbol{121} & \boldsymbol{209} & 
\boldsymbol{253} & \boldsymbol{209} & \boldsymbol{121} & \boldsymbol{44} & 
\boldsymbol{11} & 0 & 0 & 0 & 0 \\ 
12 & 54 & 166 & 329 & 463 & 461 & 331 & 164 & 56 & 10 & 2 & 0 & 0 \\ 
\boldsymbol{13} & \boldsymbol{65} & \boldsymbol{221} & \boldsymbol{494} & 
\boldsymbol{793} & \boldsymbol{923} & \boldsymbol{793} & \boldsymbol{494} & 
\boldsymbol{221} & \boldsymbol{65} & \boldsymbol{13} & 0 & 0 \\ 
\ast & \ast & \ast & \ast & \ast & \ast & \ast & \ast & \ast & \ast & \ast & 
\ast & \ast%
\end{array}%
\right] 
\]

Another property of \textbf{Alkarkhi's Triangle } to charecterize primality
of a given number .

\textbf{Theorem 2 }: \textit{p is prime if and only if all the entires of
the }$p^{th}$ \ row not exceeding the enteries of the diagonal \textit{are
(1,-1,1,-1,1.....) mod p .}

Proof : A direct consequence of the construction and \ the recurrence
formula ,

Similarly as in the previous theorem each entery of the \textit{(}$p)^{th}$
row in Alkarkhi's Triangle has the form \ 
\[
\binom{p}{d}+(-1)^{d+1} 
\]%
in $p^{th}$ row 
\[
\left[ 
\begin{array}{ccccc}
\binom{p}{1}+(-1)^{1+1}, & \binom{p}{2}+(-1)^{2+1}, & \binom{p}{2}%
+(-1)^{3+1}, & \binom{p}{d}+(-1)^{4+1}, & ...%
\end{array}%
\right] 
\]%
\[
\equiv \left[ 
\begin{array}{cccc}
0\func{mod}(p)+(-1)^{1+1}, & 0\func{mod}(p)+(-1)^{2+1}, & 0\func{mod}%
(p)+(-1)^{3+1}, & 0\func{mod}(p)+(-1)^{4+1},%
\end{array}%
\right] 
\]%
\[
\equiv \left[ 
\begin{array}{cccccc}
(-1)^{1+1}, & (-1)^{2+1}, & (-1)^{3+1}, & (-1)^{4+1}, & (-1)^{5+1}, & ...%
\end{array}%
\right] \func{mod}p 
\]%
\[
\left[ 
\begin{array}{cccccc}
1, & -1, & 1, & -1, & 1, & ...%
\end{array}%
\right] \func{mod}p 
\]%
\textbf{Alkakhi's triangle rows sum }

\textbf{Lemma 2} : The sum for the $n^{th}$ row gives power of 2 as%
\[
\sum_{d=1}^{n}K(n,d)=\left\{ 
\begin{array}{c}
2^{n}\text{\ \ \ \ \ \ \ \ \ \ \ \ \ \ }if\text{ \ \ \ n \ is odd\ \ \ \ \ \
\ \ \ \ \ \ \ \ \ } \\ 
2^{n}\text{\ }-1\text{\ \ \ \ }if\text{ \ \ \ \ \ n \ is even \ \ \ \ \ \ \
\ \ \ \ \ \ }%
\end{array}%
\right. 
\]

Proof :%
\begin{eqnarray*}
\sum_{d=1}^{n}K(n,d) &=&\sum_{d=1}^{n}\left\{ \binom{n}{d}+(-1)^{d+1}\right\}
\\
&=&\sum_{d=1}^{n}\binom{n}{d}+\sum_{d=1}^{n}(-1)^{d+1} \\
&=&(2^{n}\text{\ }-1)+\left\{ 
\begin{array}{c}
1\text{\ \ \ \ }if\text{ \ \ \ n \ is odd\ \ \ \ \ \ \ \ \ \ \ \ \ \ \ } \\ 
0\text{\ \ \ \ }if\text{ \ \ \ \ \ n \ is even \ \ \ \ \ \ \ \ \ \ \ \ \ }%
\end{array}%
\right. \text{\ } \\
\text{\ } &\text{=}&\left\{ 
\begin{array}{c}
2^{n}\text{\ \ \ \ \ \ \ \ \ \ \ \ \ \ }if\text{ \ \ \ n \ is odd\ \ \ \ \ \
\ \ \ \ \ \ \ \ \ } \\ 
2^{n}\text{\ }-1\text{\ \ \ \ }if\text{ \ \ \ \ \ n \ is even \ \ \ \ \ \ \
\ \ \ \ \ \ }%
\end{array}%
\right.
\end{eqnarray*}

Example: The sum of the 7$^{th}$ row enteries :%
\[
\begin{array}{ccccccc}
8 & 20 & 36 & 34 & 22 & 6 & 2%
\end{array}%
\]%
\ \ is given as $\ \ \ \ 8+20+36+34+22+6+2=128=2^{7}$

Example : \ The sum of the 10$^{th}$ row enteries 
\[
\begin{array}{ccccccccc}
11 & 44 & 121 & 209 & 253 & 209 & 121 & 44 & 11%
\end{array}%
\]%
is given by $11+44+121+209+253+209+121+44+11=\allowbreak 1023=2^{10}-1$%
\bigskip

\textbf{Alkakhi's Triangle rows alternating sign sums }

Lemma 3 : The alternating sign sum for the $n^{th}$ row of Alkarkhi's
Triangle is given by :%
\[
\sum\limits_{d=1}^{n}(-1)^{d+1}K(n,d)=n+1 
\]

Proof 
\begin{eqnarray*}
\sum\limits_{d=0}^{n}(-1)^{d}\binom{n}{d} &=&0 \\
&\Rightarrow &1+\sum\limits_{d=1}^{n}(-1)^{d}\binom{n}{d}=1+\sum%
\limits_{d=1}^{n}(-1)^{d}\{K(n,d)+(-1)^{d}\} \\
&=&1+\sum\limits_{d=1}^{n}(-1)^{d}\{K(n,d)+(-1)^{d}\}=1+\sum%
\limits_{d=1}^{n}\{(-1)^{d}K(n,d)+1\} \\
&=&\sum\limits_{d=1}^{n}(-1)^{d}K(n,d)+n+1=0 \\
&\Rightarrow &\sum\limits_{d=1}^{n}(-1)^{d+1}K(n,d)=n+1
\end{eqnarray*}

Example : The alternating sign sum of the \ 10$^{th}$ row enteries 
\[
\begin{array}{ccccccccccccc}
11 & 44 & 121 & 209 & 253 & 209 & 121 & 44 & 11 & 0 & 0 & 0 & 0%
\end{array}%
\]%
is given by $\ 11-44+121-209+253-209+121-44+11=10+1$

\bigskip \textbf{Alkaraji's Triangle rising diagonals and Fibonacci numbers:}

Using the well known identity 
\[
F_{n+1}=\sum\limits_{d=0}^{\left\lfloor \frac{n}{2}\right\rfloor }\binom{n-d%
}{d} 
\]

we noticed that Fibonacci numbers are located in the rising diagonals of
Alkarkhi's triangle

\[
K(n,d)=\binom{n}{d}+(-1)^{d+1} 
\]

\begin{eqnarray*}
F_{n+1} &=&1+\sum\limits_{d=1}^{\left\lfloor \frac{n}{2}\right\rfloor }%
\binom{n-d}{d} \\
&=&1+\sum\limits_{d=1}^{\left\lfloor \frac{n}{2}\right\rfloor
}(K(n-d,d)+(-1)^{d}) \\
&=&1+\sum\limits_{d=1}^{\left\lfloor \frac{n}{2}\right\rfloor
}K(n-d,d)+\left\{ 
\begin{array}{c}
-1\text{\ \ \ \ \ }if\text{ \ \ }\left\lfloor \frac{n}{2}\right\rfloor \text{
\ is odd\ \ \ \ \ \ \ \ \ \ \ \ \ \ \ } \\ 
0\text{\ \ \ \ }if\text{ \ \ }\left\lfloor \frac{n}{2}\right\rfloor \text{ \
is even \ \ \ \ \ \ \ \ \ \ \ \ \ }%
\end{array}%
\right. \\
&=&\sum\limits_{d=1}^{\left\lfloor \frac{n}{2}\right\rfloor
}K(n-d,d)+\left\{ 
\begin{array}{c}
0\text{\ \ \ \ \ }if\text{ \ }\left\lfloor \frac{n}{2}\right\rfloor \text{ \
is odd\ \ \ \ \ \ \ \ \ \ \ \ \ \ \ } \\ 
1\text{\ \ \ \ }if\text{ \ \ }\left\lfloor \frac{n}{2}\right\rfloor \text{ \
is even \ \ \ \ \ \ \ \ \ \ \ \ \ }%
\end{array}%
\right.
\end{eqnarray*}

finaly%
\[
F_{n+1}=\left\{ 
\begin{array}{c}
\sum\limits_{d=1}^{\left\lfloor \frac{n}{2}\right\rfloor }K(n-d,d)\text{\ \
\ \ \ }if\text{ \ }\left\lfloor \frac{n}{2}\right\rfloor \text{ \ is odd\ \
\ \ \ \ \ \ \ \ \ \ \ \ \ } \\ 
\sum\limits_{d=1}^{\left\lfloor \frac{n}{2}\right\rfloor }K(n-d,d)+1\text{\
\ \ \ }if\text{ \ \ }\left\lfloor \frac{n}{2}\right\rfloor \text{ \ is even
\ \ \ \ \ \ \ \ \ \ \ \ \ }%
\end{array}%
\right. 
\]

Example : To compute $F_{8}$ , observe that n =9 and $\left\lfloor \frac{7}{2%
}\right\rfloor $ = 3 \ is odd ,we apply the formula : \ 
\[
F_{8}=\sum\limits_{d=1}^{3}K(7-d,d)+1=K(6,1)+K(5,2)+K(4,3)=7+9+5=21 
\]

Or \ one can compare with rising diagonals in Alkarkhi's triangle $\ \ $%
\[
\begin{array}{cccccccc}
2 & 0 & 0 & 0 & 0 & 0 & 0 & 0 \\ 
3 & 0 & 0 & 0 & 0 & 0 & 0 & 0 \\ 
4 & 2 & 2 & 0 & 0 & 0 & 0 & 0 \\ 
5 & 5 & 5+\nearrow & 0 & 0 & 0 & 0 & 0 \\ 
6 & 9+\nearrow & 11 & 4 & 2 & 0 & 0 & 0 \\ 
7+\nearrow & 14 & 21 & 14 & \mathbf{7} & 0 & 0 & 0 \\ 
8 & 20 & 36 & 34 & 22 & 6 & 2 & 0 \\ 
9 & 27 & 57 & 69 & 57 & 27 & 9 & 0%
\end{array}%
\]

\textbf{References :}

[1] \ Rashed, Roshdi (1970--80). "Al-KarajI (or Al-Karkh\={\i}), Abu Bakr
Ibn Muhammad Ibn al Husayn". Dictionary of Scientific Biography.
NewYork:Charles Scribner's Sons. ISBN 978-0-684-10114-9.

[2] Karl Dilcher and Kenneth B.Stolasky. Apascal-Type Triangle
Characterizing Twin Primes.Amer.Math.Monthly,112(8):673-681,2005.

\end{document}